\documentclass[12pt]{article}

\usepackage[a4paper,margin=1in]{geometry}
\usepackage[T1]{fontenc}
\usepackage[utf8]{inputenc}
\usepackage{lmodern}
\usepackage{microtype}

\usepackage{amsmath,amssymb,amsthm}
\usepackage{aliascnt}
\usepackage{array,tabularx}
\usepackage{enumitem}
\usepackage[section]{placeins}
\usepackage{setspace}
\usepackage{titlesec}
\usepackage{tikz}
\usetikzlibrary{calc}
\usepackage{xcolor}
\usepackage{float}
\usepackage[
font=small,
labelfont=bf,
labelsep=period,
justification=centering
]{caption}

\usepackage[hidelinks]{hyperref}
\usepackage[nameinlink,capitalize,noabbrev]{cleveref}

\hypersetup{
	pdftitle={Turán numbers of triangular pyramids},
	pdfauthor={Anonymous Author(s)},
	pdfsubject={Turán numbers and extremal graphs for triangular pyramids},
	pdfkeywords={Turán number, triangular pyramid, extremal graph,
		lower-order term, cycle packing, join}
}

\allowdisplaybreaks
\setlist[itemize]{
	leftmargin=1.60em,
	itemsep=2pt,
	topsep=4pt
}
\setlist[enumerate]{
	leftmargin=1.95em,
	itemsep=2pt,
	topsep=4pt
}

\titleformat{\section}
{\normalfont\large\bfseries\raggedright}
{\thesection}
{0.75em}
{}
\titlespacing*{\section}
{0pt}
{3.0ex plus 1ex minus .2ex}
{1.5ex}

\makeatletter
\let\c@table\c@figure

\makeatother

\newtheorem{theorem}{Theorem}[section]

\newaliascnt{lemma}{theorem}
\newtheorem{lemma}[lemma]{Lemma}
\aliascntresetthe{lemma}

\newaliascnt{proposition}{theorem}
\newtheorem{proposition}[proposition]{Proposition}
\aliascntresetthe{proposition}

\newaliascnt{corollary}{theorem}

\aliascntresetthe{corollary}

\theoremstyle{definition}

\newaliascnt{problem}{theorem}
\newtheorem{problem}[problem]{Problem}
\aliascntresetthe{problem}

\newaliascnt{definition}{theorem}
\newtheorem{definition}[definition]{Definition}
\aliascntresetthe{definition}

\theoremstyle{remark}

\newaliascnt{remark}{theorem}

\aliascntresetthe{remark}

\crefname{theorem}{Theorem}{Theorems}
\Crefname{theorem}{Theorem}{Theorems}
\crefname{lemma}{Lemma}{Lemmas}
\Crefname{lemma}{Lemma}{Lemmas}
\crefname{proposition}{Proposition}{Propositions}
\Crefname{proposition}{Proposition}{Propositions}
\crefname{corollary}{Corollary}{Corollaries}
\Crefname{corollary}{Corollary}{Corollaries}
\crefname{problem}{Problem}{Problems}
\Crefname{problem}{Problem}{Problems}
\crefname{definition}{Definition}{Definitions}
\Crefname{definition}{Definition}{Definitions}
\crefname{remark}{Remark}{Remarks}
\Crefname{remark}{Remark}{Remarks}

\newcommand{\ex}{\operatorname{ex}}
\newcommand{\EX}{\operatorname{EX}}
\newcommand{\TP}{\mathrm{TP}}
\newcommand{\cJ}{\mathcal J}

\newcommand{\join}{\vee}
\newcommand{\eps}{\varepsilon}

\definecolor{colzero}{RGB}{150,25,25}
\definecolor{colone}{RGB}{45,105,200}
\definecolor{coltwo}{RGB}{130,190,90}

\title{Tur\'an numbers of triangular pyramids}
\author{Chuanqi Xiao$^{1,}$\footnote{Email: xiao\_chuanqi@outlook.com.},\;
	Long-Tu Yuan$^{2,3}$\footnote{Email: ltyuan@math.ecnu.edu.cn. }, \;  Xiaowen Zhang$^{3,}$\footnote{Corresponding author. Email: xiaowzhang0128@126.com.}\\
	\small $^1$ School of Science, Nanjing University of Posts and Telecommunications\\
	\small $^2$ School of Mathematical Sciences, Key Laboratory of MEA (Ministry of Education), \\  \small Shanghai Key Laboratory of PMMP,\\ 
	\small and Nantong Institute for Applied Mathematics and Artificial Intelligence \\ 
	\small $^{3}$ East China Normal University, Shanghai, 200241, China \\}

	\date{}
	
\begin{document}
\maketitle
\begin{abstract}
	Let $\TP_k$ denote the triangular pyramid with $k$ layers.  For a graph
	family $\mathcal F$, let $\ex(n,\mathcal F)$ be the maximum number of edges
	in an $n$-vertex graph containing no member of $\mathcal F$ as a subgraph. For a graph $H$, let $\cJ(H)$ be the family of subgraph-minimal graphs among all such graphs $H-I$, where $I\subseteq V(H)$ is an independent set.  For every integer $k\ge4$ and all sufficiently
	large $n$, we prove that
	\[
	\ex(n,\TP_k)
	=\max_{0\le m\le n}
	\{m(n-m)+\ex(m,\cJ(\TP_k))\}+O_k(n)
	=\left\lfloor n^2/4\right\rfloor+\Theta_k(n^{4/3}).
	\]

	Moreover, for $k=4$, we characterize all extremal
	graphs and show that 
	\[
	\ex(n,\TP_4)
	=
	\max_{0\le m\le n}
	\bigl\{m(n-m)+\ex\big(m,\cJ(\TP_4)\big)\bigr\}
	\]
	for sufficiently large $n$. 
	
	
\end{abstract}

\section{Introduction}\label{sec:introduction}

All graphs are finite and simple. 
Let $G$ be such a graph, the vertex and edge set of $G$ is denoted by $V(G)$ and $E(G)$, the number of vertices and edges in $G$ by $v(G)$ and $e(G)$, respectively. Denote by $\alpha(G)$ and $\chi(G)$ the independence number and chromatic number of $G$. Denote by $I_n$ the independent set on $n$ vertices and by $\overline{G}$ the edge complement of the graph $G$. For $x\in V(G)$ and $S\subseteq V(G)$, let $N_G(x)$ be the neighborhood of
$x$, $G[S]$ be the subgraph induced by $S$, let $d_G(x,S)=|N_G(x)\cap S|$, and write $d_G(x)=d_G(x,V(G))$,
$\delta(G)=\min_{x\in V(G)}d_G(x)$, $\Delta(G)=\max_{x\in V(G)}d_G(x)$ and $G-S=G[V(G)\setminus S]$. For a positive integer $k$, denote $[k]=\{1,2,\ldots,k\}$.

Let $A,B\subseteq V(G)$ be disjoint, and let
$e_G(A,B)$ denote the number of edges between sets $A$ and $B$ in $G$.  Two disjoint vertex sets $U$ and $W$ are completely joined in $G$ if $uw\in E(G)$ for all $u\in U$ and $w\in W$. Let $tF$ denote the disjoint union of $t$ copies of $F$. Denote by $G_1\join G_2$ the graph obtained from $G_1 \cup G_2$,  the
vertex disjoint union of graphs $G_{1}$ and $G_{2}$, and completely join $V(G_1)$ and $V(G_2)$. 

For a family $\mathcal F$ of graphs, let $\ex(n,\mathcal F)$ be
the maximum number of edges in an $n$-vertex graph containing no member of
$\mathcal F$, and let $\EX(n,\mathcal F)$ denote the corresponding family
of extremal graphs.  When $\mathcal F=\{H\}$, we write $\ex(n,H)$ and
$\EX(n,H)$.  Let $P_k$ be the path on $k$ vertices, $C_t$ be the cycle on $t$ vertices, $K_m$ be the clique on $m$ vertices.

The Erd\H{o}s--Stone--Simonovits theorem
\cite{ErdosStone1946,ErdosSimonovits1966} determines the quadratic main term
when $\chi(H)\ge3$, but the lower-order term and extremal graphs depend on the
structure of $H$.

We use the definition of triangular pyramids from
\cite{GhoshEtAl2022,ChenChenZhu2026}. 

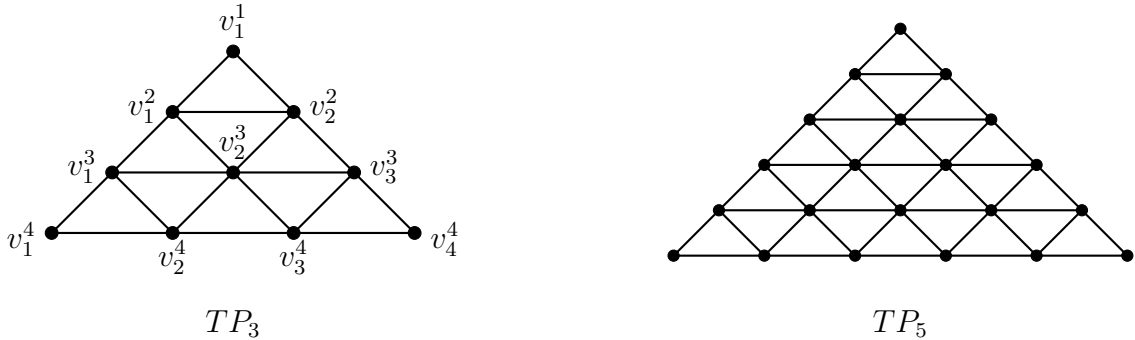
\begin{figure}[htbp]
	\centering
	\begin{tikzpicture}[scale=0.4]
		\draw[fill=black](0,0)circle(6pt);
		\draw[fill=black](-2,-2)circle(6pt);
		\draw[fill=black](2,-2)circle(6pt);
		\draw[fill=black](-4,-4)circle(6pt);
		\draw[fill=black](0,-4)circle(6pt);
		\draw[fill=black](4,-4)circle(6pt);
		\draw[fill=black](-6,-6)circle(6pt);
		\draw[fill=black](-2,-6)circle(6pt);
		\draw[fill=black](2,-6)circle(6pt);
		\draw[fill=black](6,-6)circle(6pt);
		\draw[thick](0,0)--(-6,-6)(2,-2)--(-2,-6)(4,-4)--(2,-6)(0,0)--(6,-6)(-2,-2)--(2,-6)(-4,-4)--(-2,-6);
		\draw[thick](-2,-2)--(2,-2)(-4,-4)--(4,-4)(-6,-6)--(6,-6);
		\node at (0,-9) {$TP_3$};
		\node at (0,1) {$v_1^{1}$};
		\node at (-3,-1.8) {$v_1^{2}$};
		\node at (3,-1.8) {$v_2^{2}$};
		\node at (-5,-3.8) {$v_1^{3}$};
		\node at (0,-3) {$v_2^{3}$};
		\node at (5,-3.8) {$v_3^{3}$};
		\node at (-7,-6.2) {$v_1^{4}$};
		\node at (-2,-6.8) {$v_2^{4}$};
		\node at (2,-6.8) {$v_3^{4}$};
		\node at (7,-6.2) {$v_4^{4}$};
		
	\end{tikzpicture} \qquad\qquad\qquad
	\begin{tikzpicture}[scale=0.3]
		\draw[fill=black](0,0)circle(7pt);
		\draw[fill=black](-2,-2)circle(7pt);
		\draw[fill=black](2,-2)circle(7pt);
		\draw[fill=black](-4,-4)circle(7pt);
		\draw[fill=black](0,-4)circle(7pt);
		\draw[fill=black](4,-4)circle(7pt);
		\draw[fill=black](-6,-6)circle(7pt);
		\draw[fill=black](-2,-6)circle(7pt);
		\draw[fill=black](2,-6)circle(7pt);
		\draw[fill=black](6,-6)circle(7pt);
		\draw[fill=black](-8,-8)circle(7pt);
		\draw[fill=black](-4,-8)circle(7pt);
		\draw[fill=black](0,-8)circle(7pt);
		\draw[fill=black](4,-8)circle(7pt);
		\draw[fill=black](8,-8)circle(7pt);
		\draw[fill=black](-10,-10)circle(7pt);
		\draw[fill=black](-6,-10)circle(7pt);
		\draw[fill=black](-2,-10)circle(7pt);
		\draw[fill=black](2,-10)circle(7pt);
		\draw[fill=black](6,-10)circle(7pt);
		\draw[fill=black](10,-10)circle(7pt);
		\draw[thick](0,0)--(-10,-10)(2,-2)--(-6,-10)(4,-4)--(-2,-10)(6,-6)--(2,-10)(8,-8)--(6,-10);
		\draw[thick](0,0)--(10,-10)(-2,-2)--(6,-10)(-4,-4)--(2,-10)(-6,-6)--(-2,-10)(-8,-8)--(-6,-10);
		\draw[thick](-2,-2)--(2,-2)(-4,-4)--(4,-4)(-6,-6)--(6,-6)(-8,-8)--(8,-8)(-10,-10)--(10,-10);
		\node at (0,-13) {$TP_5$};
	\end{tikzpicture}
	\caption{Triangular Pyramids with $3$ and $5$ layers respectively.}
	\label{PT3PT5}
\end{figure}

\begin{definition}\label{TPk}
	The Triangular Pyramid with $k$ layers, denoted by $TP_k$, is defined as follows: Draw $k+1$ paths in layers such that the $i$-th layer's path is an $i$-vertex path $P_i=v_1^iv_2^i\cdots v_i^i$, where $i\in [k+1]$. Then, for any consecutive layers, say the $i$-th and $(i+1)$-th, $v_r^iv_r^{i+1}\in E(TP_k)$ and  $v_r^iv_{r+1}^{i+1}\in E(TP_k)$, $r\in [i]$ (see Figure~\ref{PT3PT5}). In particular, $\TP_1=K_3$.
\end{definition}

Liu~\cite{Liu2013} determined $\ex(n,\TP_2)$ and the extremal graphs for all sufficiently large $n$, and Xiao, Katona, Xiao and Zamora
\cite{XiaoKatonaXiaoZamora2022} subsequently solved the problem for every
$n$.  Ghosh, Gy\H{o}ri, Paulos, Xiao and Zamora~\cite{GhoshEtAl2022} proved
$\ex(n,\TP_3)=n^2/4+n+o(n)$ and conjectured that
$\ex(n,\TP_4)=n^2/4+\Theta(n^{4/3})$.  Chen, Chen and
Zhu~\cite{ChenChenZhu2026} proved
$\ex(n,\TP_4)=n^2/4+\Theta(n^{4/3})$.
For $k=4$, we obtain an exact join formula and determine all extremal graphs of
sufficiently large order.

We use structural results of Simonovits
\cite{Simonovits1968,Simonovits1974,Simonovits1983}, in particular the form of the Erd\H{o}s--Simonovits structure theorem stated in
\cite{GriggsSimonovitsThomas1998}.  Our approach is also related to the framework for lower-order Tur\'an terms developed by Ma and
Qiu~\cite{MaQiu2020}.  Peng, Song and Yuan~\cite{PengSongYuan2026} studied
when such extremal problems admit exact join formulas, while Xu~\cite{Xu2026}
showed that a superlinear lower-order term alone does not force an exact join
formula.  We also use the $O(n^{4/3})$ extremal estimate for honeycomb
graphs due to Gao, Janzer, Liu and Xu~\cite{GaoJanzerLiuXu2023}.

Given a graph $H$, for each independent set $I$ of $H$, consider the graph $H-I$. Let $\cJ(H)$ be the family of subgraph-minimal graphs among all such graphs $H-I$.

The lower bound in the following theorem follows from
\Cref{prop:join-criterion-general}. We prove that $\ex(n,\TP_k)$
exceeds this bound by at most $O_k(n)$ and describe the extremal graphs up to $O_k(n)$ edge modifications.

\begin{theorem}\label{thm:structure-intro}
	For every integer $k\ge4$, there exist a constant $C_k>0$ and an
	integer $n_k$ such that, for every $n\ge n_k$,
	\[
	0\le \ex(n,\TP_k)
	-\max_{0\le m\le n}\{m(n-m)+\ex(m,\cJ(\TP_k))\}\le C_kn.
	\]
	Moreover, every $G\in\EX(n,\TP_k)$ differs in at most $C_kn$ edges from a graph $R\join I_{n-m}$, where $R$ is an
	$m$-vertex $\cJ(\TP_k)$-free graph and $|m-n/2|\le C_kn^{1/2}$.
\end{theorem}

In particular, for every $k\ge4$,
$\ex(n,\TP_k)=\left\lfloor n^2/4\right\rfloor
+\Theta_k(n^{4/3}).$ When $k=4$, the linear error disappears: for sufficiently large $n$,
every graph in $\EX(n,\TP_4)$ is a join of a $\cJ(\TP_4)$-free graph and an independent set.

\begin{theorem}\label{exactp4}
	For sufficiently large $n$,
	\[
	\ex(n,\TP_4)=\max_{0\le m\le n}
	\{m(n-m)+\ex(m,\cJ(\TP_4))\}.
	\]
	Moreover, $\EX(n,\TP_4)$ consists precisely of the graphs $R\join I_{n-m}$, where $R\in\EX(m,\cJ(\TP_4))$ and the maximum above is attained at $m$.
\end{theorem}

The remainder of the paper is organized as follows.
In \Cref{sec:preliminaries}, we establish the auxiliary results for
$\TP_k$, prove a general structural theorem, and apply it to $\TP_k$.
In \Cref{sec:exact-four-layer}, we prove \Cref{exactp4} and characterize
the extremal $\TP_4$-free graphs.
Finally, \Cref{sec:conclusion} presents the remaining open problems.
The determination of $\cJ(\TP_4)$ is the only computer-assisted step;
details are given in \Cref{app:verification}.

\section{Preliminaries and auxiliary estimates}\label{sec:preliminaries}

\makeatletter
\tagsleft@false
\let\veqno\@@eqno
\makeatother

\titleformat{\subsection}[block]
{\normalfont\normalsize\bfseries\raggedright}{\thesubsection}{0.65em}{}
\titlespacing*{\subsection}{0pt}{2.3ex plus .8ex minus .2ex}{.8ex}
\titleformat{\subsubsection}[block]
{\normalfont\normalsize\itshape\raggedright}{\thesubsubsection}{0.65em}{}
\titlespacing*{\subsubsection}{0pt}{1.8ex plus .6ex minus .2ex}{.6ex}

We first use the proper $3$-coloring of $\TP_k$ to determine the bipartite members of $\cJ(\TP_k)$. We then determine the maximum number of pairwise vertex-disjoint copies of $C_6$ in the graphs $B_{k,r}$ and prove that, for every $k\ge4$, $\ex(n,\cJ(\TP_k))=\Theta_k(n^{4/3})$. Finally, we estimate the integers $m$ maximizing $m(n-m)+\ex(m,\cJ(\TP_k))$. These results will be used in the proof of \Cref{thm:structure-intro}.

By the definition of $TP_k$, we have $v(\TP_k)=\binom{k+2}{2}$ and $e(\TP_k)=3k(k+1)/2$. We color $v_j^i$ by the residue of $i+j$ modulo $3$.
The values of $i+j$ at the endpoints of the three types of
edges in $\TP_k$ differ by $1$, $1$, and $2$, respectively.
Thus this is a proper $3$-coloring. Since $\TP_k$ contains a triangle, $\chi(\TP_k)=3$.  The coloring is
unique up to relabeling: once one triangular face is colored, the color on
every adjacent face is forced, and the triangular faces form a connected
region.

For $r\in\{0,1,2\}$, let $V_r=\{v_j^i\in V(\TP_k):i+j\equiv r\pmod3\}$
and let $B_{k,r}=\TP_k-V_r$ be the two-color induced subgraph obtained by
deleting the color class $V_r$ from $\TP_k$.  The two remaining color
classes form a bipartition of $B_{k,r}$.

These graphs identify the bipartite members of
$\cJ(\TP_k)$. If $\TP_k-I\in\cJ(\TP_k)$ is bipartite, then its
two parts together with the independent set $I$ form a proper
$3$-coloring of $\TP_k$.  The uniqueness of the coloring shows that
$\TP_k-I$ is one of $B_{k,0},B_{k,1},B_{k,2}$.  Conversely, $B_{k,r}$ is obtained from $\TP_k$ by deleting
the independent set $V_r$.
By definition, at least one graph in $\cJ(\TP_k)$
is a subgraph of $B_{k,r}$.

The family $\cJ(H)$ also determines which graphs
can be joined to an independent set without creating $H$.  The following
proposition makes this correspondence precise and gives the lower bound
used throughout the paper.

\begin{proposition}\label{prop:join-criterion-general}
	Let $H$ be a graph with $\chi(H)\ge3$, let $R$ be a graph, and let
	$b$ be a nonnegative integer.
	\begin{enumerate}[label=\textup{(\roman*)}]
		\item If $R$ is $\cJ(H)$-free, then $R\join I_b$ is $H$-free.
		\item If $b\ge\alpha(H)$, then $ R\join I_b$ is $H$-free if and only if $R$ is  $\cJ(H)$-free.
	\end{enumerate}
	Consequently, for every $n$,
	$\ex(n,H)\ge
	\max_{0\le m\le n}\{m(n-m)+\ex(m,\cJ(H))\}.$
\end{proposition}

\begin{proof}
	Suppose that $H\subseteq R\join I_b$. Let $I=V(H)\cap V(I_b)$, then  $I$ is an independent set of $H$. Moreover,
	$H-I\subseteq R$. By the definition of $\cJ(H)$, the graph $H-I$
	contains a member of $\cJ(H)$. Therefore, $R$ contains a member of
	$\cJ(H)$, this proves \textup{(i)}.
	
	Now, suppose that $R$ contains a graph $M\in\cJ(H)$. By the definition of $\cJ(H)$, there is an independent set $I\subseteq V(H)$ such that $M\cong H-I$. As $|I|\le\alpha(H)\le b$, choose $|I|$ vertices from $I_b$. Since every vertex of $R$ is adjacent to every
	vertex of $I_b$, these vertices together with $M$ contains a copy of $H$. Combining with \textup{(i)}, this proves \textup{(ii)}.
	
	Finally, since $\chi(H)\ge3$, the graph $H-I$ contains an edge for every
	independent set $I\subseteq V(H)$. Otherwise, both $I$ and
	$V(H)\setminus I$ would be independent, and $H$ would be bipartite.
	Thus, for every $m\leq n$, there exists an $m$-vertex $\cJ(H)$-free graph $R$ with
	$\ex(m,\cJ(H))$ edges. By \textup{(i)}, $R\join I_{n-m}$ is an $n$-vertex
	$H$-free graph with $m(n-m)+\ex(m,\cJ(H))$ edges.
	Since this construction is valid for every $0\le m\le n$,
	we have
	$ \ex(n,H)\ge
	\max_{0\le m\le n}
	\bigl\{m(n-m)+\ex(m,\cJ(H))\bigr\}.$
\end{proof}

Let $\nu_6(G)$ denote the maximum
number of pairwise vertex-disjoint copies of $C_6$ in $G$.  We now
determine $\nu_6(B_{k,r})$. Next, the fact that this
number is positive for every $r$ when $k\ge4$ will give the lower bound for
$\ex(n,\cJ(\TP_k))$.

\begin{proposition}\label{prop:c6-values}
	Let $k\ge1$.
	\begin{enumerate}[label=\textup{(\roman*)}]
		\item If $k=3q+1$ or $k=3q+2$, where $q\ge0$, then
		\[
		\nu_6(B_{k,0})=\nu_6(B_{k,1})=\nu_6(B_{k,2})
		=\binom{q+1}{2}.
		\]
		\item If $k=3t$, where $t\ge1$, then
		\[
		\nu_6(B_{k,0})=\nu_6(B_{k,1})=\binom t2,
		\qquad
		\nu_6(B_{k,2})=\binom{t+1}{2}.
		\]
	\end{enumerate}
	Consequently,
	\[
	\min_{r\in\{0,1,2\}}\nu_6(B_{k,r})
	=\binom{\lfloor(k-1)/3\rfloor+1}{2}.
	\]
\end{proposition}

\begin{proof}
	Let $r\in\{0,1,2\}$.  Extend the notation
	$v_j^i$ to all integer pairs $(i,j)$ and retain the three edge rules in the definition of $\TP_k$.  Color $v_j^i$ by $i+j$ modulo $3$ and delete all vertices of color $r$.  Every bounded face of the resulting plane graph is a $6$-cycle surrounding one deleted vertex.
	
	Every $6$-cycle in this plane graph is the boundary of one such face.
	Let $C$ be a $6$-cycle, and let $P$ be the plane subgraph formed by
	the boundaries of all faces in the closed disk bounded by $C$.  These faces
	form a connected region whose outer boundary is $C$.  If there are $f$ such
	faces, then Euler's formula and the sum of the face lengths give
	$v(P)-e(P)+f=1$, and $2e(P)=6f+6$.
	Thus $e(P)=3f+3$ and $v(P)=2f+4$.  The $2f-2$ vertices not on $C$ have degree three in $P$, so the total degree in $P$ of the six vertices on $C$ is $2e(P)-3(2f-2)=12$. Each of these six vertices already has degree two on $C$.  Hence none is
	joined to a vertex inside $C$.  If $f>1$, then $P$ has
	$v(P)-6=2f-2>0$ vertices inside $C$, contradicting the connectedness of $P$, hence $f=1$.
	
	The face surrounding $v_j^i$ has boundary $ v_{j-1}^{i-1},\ v_j^{i-1},\ v_{j+1}^{i},\
	v_{j+1}^{i+1},\ v_j^{i+1},\ v_{j-1}^{i}.$
	All six vertices lie in $\TP_k$ exactly when $2\le i\le k$ and
	$2\le j\le i-1$.  Its boundary lies in $B_{k,r}$ exactly when
	$i+j\equiv r\pmod3$.  Thus the copies of $C_6$ in $B_{k,r}$ correspond
	precisely to the vertices $v_j^i$ satisfying these three conditions.
	
	For each such copy, let
	$a=j-1$, $b=i-j$ and $c=k+1-i$.
	Thus, the copies of $C_6$ in $B_{k,r}$ are represented by the triples in $\Lambda_{k,r}$. We define
	\[
	\Lambda_{k,r}
	=\{(a,b,c):a,b,c\text{ are positive integers},\
	a+b+c=k,\ k+a-c+2\equiv r\pmod3\}.
	\]
	Let $p=(a,b,c)$ and $p'=(a',b',c')$ correspond to two distinct copies
	centered at $v_j^i$ and $v_{j'}^{i'}$, respectively. If the center has coordinates $(i,j)$, then the six vertices of the cycle have coordinates
	$ (-1,-1),\ (-1,0),\ (0,1),\ (1,1),\ (1,0),\ (0,-1)$. 
	
	Moreover, $(i'-i)+(j'-j)\equiv0\pmod3$ because the two centers have the
	same color. Comparing their vertex sets, we see that two copies of $C_6$ intersect if and only if
	$(i'-i,j'-j)\in
	\{\pm(1,2),\pm(2,1),\pm(1,-1)\}$ holds.
	Moreover,
	\[
	p'-p=(a'-a,b'-b,c'-c)
	=(j'-j,(i'-j')-(i-j),i-i').
	\]
	Therefore, two distinct triples $p,p'\in\Lambda_{k,r}$ correspond to
	intersecting copies of $C_6$ if and only if
	\begin{equation}\label{eq:overlap-vectors}
		p'-p\in
		\{\pm(2,-1,-1),\ \pm(-1,2,-1),\ \pm(-1,-1,2)\}.
	\end{equation}
	
	By the correspondence above and \eqref{eq:overlap-vectors}, a collection of
	copies of $C_6$ in $B_{k,r}$ is pairwise vertex-disjoint if and only if the
	difference of any two corresponding triples in $\Lambda_{k,r}$ is not one
	of the vectors listed in \eqref{eq:overlap-vectors}. Hence,
	$\nu_6(B_{k,r})$ is the maximum size of such a collection. To determine
	this maximum, let $\sigma\in\{1,2,3\}$ and let
	$k=3q+\sigma$, where
	$q=\lfloor(k-1)/3\rfloor$.
	Here, $\sigma=3$ if and only if $3\mid k$. If $q=0$, then
	$k\in\{1,2,3\}$.  If $k=1$ or $2$, there are no positive integers $a,b,c$ satisfying $a+b+c=k$. If $k=3$, the only solution is $(a,b,c)=(1,1,1)$, which lies in
	$\Lambda_{3,2}$.  Thus the proposition holds when $q=0$, and we may assume
	$q\ge1$.
	For nonnegative integers $x,y,z$ with $x+y+z=q-1$, define $A_{x,y,z}=(3x+1,3y+2,3z+\sigma)$.
	These triples belong to $\Lambda_{k,0}$. Their pairwise differences have
	all three coordinates divisible by $3$, so none occurs in
	\eqref{eq:overlap-vectors}.  They give $\binom{q+1}{2}$ pairwise
	vertex-disjoint copies of $C_6$ in $B_{k,0}$.

	To prove the upper bound, let $u=(2,-1,-1)$ and $w=(1,1,-2)$.
	For the same triples $A_{x,y,z}$, consider the sets $\Lambda_{k,0}\cap
	\{A_{x,y,z},A_{x,y,z}+u,A_{x,y,z}+w\}$.
	Their three possible elements have pairwise differences among the vectors in
	\eqref{eq:overlap-vectors}, so at most one of them can be selected.
	
	We show that every triple in $\Lambda_{k,0}$ belongs to exactly one of
	these sets. Let $p=(a,b,c)\in\Lambda_{k,0}$. The defining congruences give
	\[
	(a,b,c)\equiv
	\begin{cases}
		(1,2,\sigma),&a\equiv1,\\
		(0,1,\sigma+2),&a\equiv0,\\
		(2,0,\sigma+1),&a\equiv2
	\end{cases}
	\pmod3.
	\]
	
	In the three cases, the unique triple of the form
	$A_{x,y,z}$ is $p$, $p-u$, or $p-w$.  In the second case $a\ge3$, so
	$p-u$ has positive coordinates.  In the third case $a\ge2$ and $b\ge3$,
	so $p-w$ has positive coordinates.  There are $\binom{q+1}{2}$ choices of
	$(x,y,z)$, and hence
	$\nu_6(B_{k,0})=\binom{q+1}{2}$.
	
	The coordinate permutations
	$(a,b,c)\mapsto(b,c,a)$ and $(a,b,c)\mapsto(c,b,a)$ preserve the vectors
	in \eqref{eq:overlap-vectors}. If $(a,b,c)\in\Lambda_{k,r}$, then $k+b-a+2\equiv r+k\pmod3$. Hence $(a,b,c)\mapsto(b,c,a)$ is a bijection from $\Lambda_{k,r}$ to
	$\Lambda_{k,s}$, where $s\equiv r+k\pmod3$.
	If $3\nmid k$, this map cyclically permutes the three residues and preserves
	whether the corresponding copies of $C_6$ intersect.  Hence the three values
	of $\nu_6(B_{k,r})$ are equal.  This proves part~\textup{(i)}.
	
	Let $k=3t$, then $q=t-1$ and $\sigma=3$, so the preceding argument gives
	$\nu_6(B_{3t,0})=\binom t2$.  If $(a,b,c)\in\Lambda_{3t,r}$, then $3t+c-a+2\equiv1-r\pmod3$. Hence $(a,b,c)\mapsto(c,b,a)$ is a bijection from
	$\Lambda_{3t,r}$ to $\Lambda_{3t,s}$, where
	$s\equiv1-r\pmod3$.
	It preserves the set of differences in \eqref{eq:overlap-vectors}, so
	$\nu_6(B_{3t,1})=\binom t2$.
	
	It remains to determine $\nu_6(B_{3t,2})$.  A triple
	$(a,b,c)\in\Lambda_{3t,2}$ has $a\equiv c\pmod3$; since
	$a+b+c\equiv0\pmod3$, all three coordinates have the same residue.
	For nonnegative integers $x,y,z$ satisfying $x+y+z=t-1$, let
	$A'_{x,y,z}=(3x+1,3y+1,3z+1)$.
	These $\binom{t+1}{2}$ triples give pairwise vertex-disjoint copies of
	$C_6$.  For each such triple $A'_{x,y,z}$, consider
	$\Lambda_{3t,2}\cap
	\{A'_{x,y,z},A'_{x,y,z}+u,A'_{x,y,z}+w\}$.
	These sets partition $\Lambda_{3t,2}$.  For a triple $p=(a,b,c)$ whose
	common residue is $1$, $0$, or $2$, respectively, the unique triple of the
	form $A'_{x,y,z}$ is $p$, $p-u$, or $p-w$.  In the second case all
	coordinates of $p$ are at least $3$, so $p-u$ has positive coordinates; in
	the third they are at least $2$, so $p-w$ has positive coordinates.  Each
	set contains at most one triple from a collection
	corresponding to vertex-disjoint copies of $C_6$.  Hence $\nu_6(B_{3t,2})=\binom{t+1}{2}$ holds,
	which proves \textup{(ii)} for all $t\ge1$.  Taking the minimum of the
	three values completes the proof.
\end{proof}

We next prove that
$\ex(n,\cJ(\TP_k))=\Theta_k(n^{4/3})$.  The lower bound follows from
the copies of $C_6$ found above and the identification of the bipartite
members of $\cJ(\TP_k)$.  For the upper bound, we apply a theorem of Gao, Janzer, Liu and Xu~\cite{GaoJanzerLiuXu2023} to a graph containing
$B_{k,r}$ as a subgraph.

We use the classical estimate
$\ex(n,C_6)=\Theta(n^{4/3})$~\cite{BondySimonovits1974,FurediNaorVerstraete2006}.  For every positive
integer $t$ and every $n\ge1$, we show
\[
\ex(n,C_6)\le \ex(n,tC_6)
\le \ex(n,C_6)+6(t-1)n.
\]
The lower bound is trivial.  For the upper
bound, take a maximal family of vertex-disjoint copies of $C_6$ in a
$tC_6$-free graph $G$.  Their union $U$ has at most $6(t-1)$ vertices and
$G-U$ is $C_6$-free, hence we get
$ e(G)\le e(G-U)+|U|n
\le \ex(n,C_6)+6(t-1)n$.
Consequently,
$\ex(n,tC_6)=\Theta_t(n^{4/3})$.

We use the finite honeycomb graphs introduced by Gao, Janzer, Liu and
Xu~\cite{GaoJanzerLiuXu2023}. We define the finite honeycomb graph $H_{s,\ell}$ as follows. For an odd integer $s\geq 1$ and even $\ell\geq 2$, let $H_{s,\ell}$ be the graph with vertex set $V(H_{s,\ell})=\{x_{i,j}:1\leq i\leq s, 1\leq j\leq \ell\}$, where $x_{s,1}=x_{s,3}=\cdots=x_{s,\ell-1}=u$ and $x_{1,2}=x_{1,4}=\cdots=x_{1,\ell}=v$ (but all the other vertices are distinct) and edge set
\[
\begin{aligned}
	E(H_{s,\ell})={}&
	\{x_{i,j}x_{i,j+1}:1\le i\le s,\ 1\le j\le \ell-1\}\\
	&\cup
	\{x_{2i-1,j}x_{2i,j}:1\le i\le (s-1)/2,\
	1\le j\le \ell,\ j\text{ is odd}\}\\
	&\cup
	\{x_{2i,j}x_{2i+1,j}:1\le i\le (s-1)/2,\
	1\le j\le \ell,\ j\text{ is even}\}.
\end{aligned}
\]

\begin{figure}[htbp]
	\centering
	\definecolor{pathgreen}{RGB}{24,130,82}
	\definecolor{oddblue}{RGB}{36,92,170}
	\definecolor{evenred}{RGB}{190,54,48}
	\begin{tikzpicture}[
		scale=0.85,
		transform shape,
		x=1cm,
		y=1cm,
		path edge/.style={draw=pathgreen, line width=1.05pt},
		odd edge/.style={draw=oddblue, line width=1.05pt},
		even edge/.style={draw=evenred, line width=1.05pt},
		vertex/.style={circle, fill=black, inner sep=1.35pt},
		boundary vertex/.style={circle, fill=black, inner sep=1.65pt},
		vertex label/.style={font=\small}
		]
		
		\foreach \j in {1,3,...,11}{
			\foreach \i in {1,...,6}{
				\coordinate (x-\i-\j) at ({1.18*(\i-1)},{0.58*(\j-1)});
			}
		}
		\foreach \j in {2,4,...,12}{
			\foreach \i in {2,...,7}{
				\coordinate (x-\i-\j) at ({1.18*(\i-1.5)},{0.58*(\j-1)});
			}
		}
		
		\coordinate (v) at (-1.08,{0.58*5.5});
		\coordinate (u) at ({1.18*6.0},{0.58*5.5});

		\foreach \i in {2,...,6}{
			\foreach \j in {1,...,11}{
				\pgfmathtruncatemacro{\jp}{\j+1}
				\draw[path edge] (x-\i-\j)--(x-\i-\jp);
			}
		}
		\foreach \j in {1,3,...,11}{
			\draw[path edge] (v)--(x-1-\j);
		}
		\foreach \j in {2,4,...,12}{
			\draw[path edge] (u)--(x-7-\j);
		}

		\foreach \j in {1,3,...,11}{
			\foreach \a/\b in {1/2,3/4,5/6}{
				\draw[odd edge] (x-\a-\j)--(x-\b-\j);
			}
		}

		\foreach \j in {2,4,...,12}{
			\foreach \a/\b in {2/3,4/5,6/7}{
				\draw[even edge] (x-\a-\j)--(x-\b-\j);
			}
		}
		
		\foreach \j in {1,3,...,11}{
			\foreach \i in {1,...,6}{\node[vertex] at (x-\i-\j) {};}
		}
		\foreach \j in {2,4,...,12}{
			\foreach \i in {2,...,7}{\node[vertex] at (x-\i-\j) {};}
		}
		\node[boundary vertex] at (v) {};
		\node[boundary vertex] at (u) {};
		
		\node[vertex label, left=3pt] at (v) {$v$};
		\node[vertex label, right=3pt] at (u) {$u$};
		\foreach \i in {1,...,6}{
			\node[vertex label, below=4pt] at (x-\i-1) {$x_{\i,1}$};
		}
		\foreach \i in {2,...,7}{
			\node[vertex label, above=4pt] at (x-\i-12) {$x_{\i,12}$};
		}
	\end{tikzpicture}
	\caption{The graph $H_{7,12}$.}
	\label{fig:honeycomb-H712}
\end{figure}
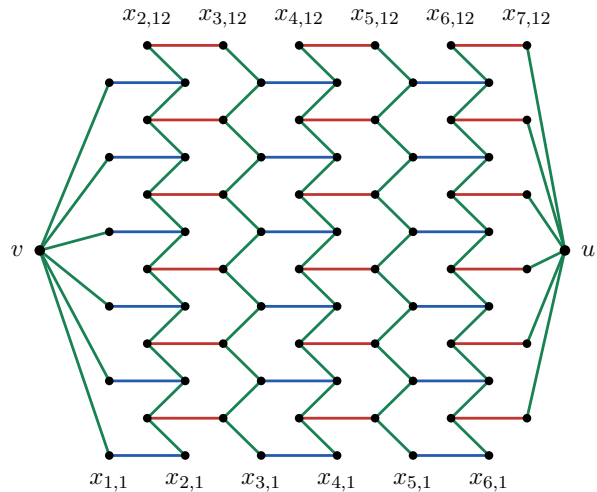

In Figure~\ref{fig:honeycomb-H712}, the edges
$x_{i,j}x_{i,j+1}$ are green. For odd $j$, the edges
$x_{2i-1,j}x_{2i,j}$ are blue; for even $j$, the edges
$x_{2i,j}x_{2i+1,j}$ are red.

Their result implies that, for every $k\ge4$,
\begin{equation}\label{thm:honeycomb-extremal}
	\ex(n,H_{2k+9,\,2k+8})=O_k(n^{4/3}).
\end{equation}
\

We now determine the order of $\ex(n,\cJ(\TP_k))$. The lower bound
follows from extremal constructions avoiding disjoint copies of $C_6$.
For the upper bound, we show that $B_{k,r}$ is a subgraph of a suitable
finite honeycomb graph and then apply
\Cref{thm:honeycomb-extremal}.

\begin{proposition}\label{prop:join-family-comparison}
	For every $k\ge4$,
	$\ex(n,\cJ(\TP_k))=\Theta_k(n^{4/3}).$
\end{proposition}

\begin{proof}
	Let
	\[
	t=\binom{\lfloor(k-1)/3\rfloor+1}{2}.
	\]
	Then $t$ is positive for $k\ge4$.  Let
	$Q\in\EX(n,tC_6)$, and let $Q'$ be a spanning bipartite subgraph obtained
	from a maximum cut of $Q$.  Then $e(Q')\ge e(Q)/2$, and $Q'$ remains
	$tC_6$-free.  Every bipartite member of $\cJ(\TP_k)$ is isomorphic to some $B_{k,r}$ and
	therefore contains $tC_6$ by \Cref{prop:c6-values}. Since $Q'$ is bipartite
	and $tC_6$-free, it is $\cJ(\TP_k)$-free. Moreover, we have
	$ \ex(n,\cJ(\TP_k))
	\ge\frac12\ex(n,tC_6)
	=\Omega_k(n^{4/3})$.
	
	For the upper bound, fix $r\in\{0,1,2\}$.  We show that $B_{k,r}$ is isomorphic to a subgraph of $H_{2k+9,2k+8}$. For each $v_j^i\in V(B_{k,r})$, denote its corresponding vertex in
	$H_{2k+9,\,2k+8}$ by $x_{\alpha,\beta}$, where
	$\beta=i+k+5-r$
	and
	\[
	\alpha=
	\begin{cases}
		k+4+(i-2j+2-r)/3,
		& i+j-r\equiv1\pmod3,\\
		k+4+(i-2j+1-r)/3,
		& i+j-r\equiv2\pmod3.
	\end{cases}
	\]
	Since $v_j^i\notin V_r$, we have
	$i+j-r\equiv1$ or $2\pmod3$. Thus exactly one of the two cases applies,
	and $\alpha$ is an integer. Moreover, $\beta$ uniquely determines $i$, and once $i$ is known,
	$\alpha$ uniquely determines $j$. Therefore, distinct vertices of $B_{k,r}$
	correspond to distinct vertices of $H_{2k+9,\,2k+8}$.
	
	It follows directly from the definitions of $\TP_k$ and $H_{s,\ell}$
	that adjacent vertices of $B_{k,r}$ correspond to adjacent vertices of
	$H_{2k+9,\,2k+8}$. Since $1\le j\le i\le k+1$, the corresponding
	coordinates satisfy
	$2\le\alpha\le2k+8$ and $2\le\beta\le2k+7$. Thus none of these
	vertices lies in either of the two boundary columns where vertex
	identifications occur. Consequently, $B_{k,r}$ is isomorphic to a
	subgraph of $H_{2k+9,\,2k+8}$.
	
	Since $2k+9$ is odd and $2k+8$ is even, by \eqref{thm:honeycomb-extremal}, we get
	\[
	\ex(n,\cJ(\TP_k))
	\le \ex(n,B_{k,r})
	\le \ex(n,H_{2k+9,\,2k+8})
	=O_k(n^{4/3}).
	\]\end{proof}

Let $m$ be an integer such that $m(n-m)+\ex(m,\cJ(\TP_k))
=
\max_{0\le s\le n}
\bigl\{s(n-s)+\ex(s,\cJ(\TP_k))\bigr\}$. We next show that every such $m$ is close to $n/2$. This estimate will
be used to refine the bound for $\ex(n,\TP_k)$ and to describe the
extremal graphs.

We first establish a property of the graphs in $\cJ(\TP_k)$ that will be used below. Let $I\subseteq V(\TP_k)$ be an independent set. We claim
that $\TP_k-I$ is connected and has no isolated vertex. The
neighborhood of a boundary vertex induces a path, and that of an interior
vertex induces a $6$-cycle.  If a path between two vertices outside $I$
passes through $x\in I$, replace the part through $x$ by a path between its
two neighboring vertices inside $N_{\TP_k}(x)$.  Repeating this operation
gives a path in $\TP_k-I$.  Finally, every vertex lies in a triangle, so no
vertex outside $I$ becomes isolated after deleting $I$.

It follows that every member of $\cJ(\TP_k)$ is connected and has no
isolated vertex.  Taking disjoint unions of extremal graphs therefore gives
\[
\ex(a+b,\cJ(\TP_k))
\ge\ex(a,\cJ(\TP_k))+\ex(b,\cJ(\TP_k))
,~\text{where}~a,b\ge0.
\]
Since every member of $\cJ(\TP_k)$ has no isolated vertex, adding an
isolated vertex to a $\cJ(\TP_k)$-free graph does not create a member of
$\cJ(\TP_k)$. Hence $\ex(m+1,\cJ(\TP_k))\ge \ex(m,\cJ(\TP_k))$ holds.

Now let $R\in\EX(m,\cJ(\TP_k))$. Then $R$ has a vertex $v$ with $ d_R(v)\le 2e(R)/m$. Moreover, the graph $R-v$ is still $\cJ(\TP_k)$-free. Hence,
\[
\ex(m-1,\cJ(\TP_k))
\ge e(R-v)
\ge \ex(m,\cJ(\TP_k))
-\frac{2\ex(m,\cJ(\TP_k))}{m}.
\]
Therefore,
\begin{equation}\label{eq:increment-bound}
	0\le\ex(m,\cJ(\TP_k))-\ex(m-1,\cJ(\TP_k))
	\le\frac{2\ex(m,\cJ(\TP_k))}{m}
	\qquad(m\ge1).
\end{equation}
Combining with \Cref{prop:join-family-comparison}, we get
$
0\le
\ex(m,\cJ(\TP_k))-\ex(m-1,\cJ(\TP_k))
\le C_km^{1/3}$
for some constant $C_k>0$.

We now use the preceding estimate to show that every maximizing integer
$m$ is close to $n/2$. We then prove a similar conclusion when
$m(n-m)+\ex(m,\cJ(\TP_k))$ is at most $O_k(n)$ below its maximum.

\begin{lemma}\label{lem:join-maximizers}
	Fix an integer $k\ge4$. There exists a constant $C_k>0$ such that,
	for all sufficiently large $n$, every integer
	$m\in\{0,\ldots,n\}$ maximizing $m(n-m)+\ex(m,\cJ(\TP_k))$ satisfies
	\[
	\left\lfloor\frac n2\right\rfloor
	\le m
	\le \frac n2+C_kn^{1/3}.
	\]
	
	Moreover, for every $D>0$, there exists a constant $C_{k,D}>0$ such
	that, for all sufficiently large $n$, every integer
	$a\in\{0,\ldots,n\}$ satisfying $ a(n-a)+\ex(a,\cJ(\TP_k))
	\ge{}
	\max_{0\le s\le n}
	\bigl\{s(n-s)+\ex(s,\cJ(\TP_k))\bigr\}-Dn$ satisfies
	$ \left|a-\frac n2\right|\le C_{k,D}n^{1/2}.
	$
\end{lemma}

\begin{proof}
	Let $m$ attain the maximum. Since
	$\ex(s,\cJ(\TP_k))$ is nondecreasing, we have
	$m\ge\lfloor n/2\rfloor$. Maximality at $m$ and
	\eqref{eq:increment-bound} give
	\[
	2m-n-1
	\le \ex(m,\cJ(\TP_k))-\ex(m-1,\cJ(\TP_k))
	=O_k(n^{1/3}),
	\]
	which proves the upper bound.
	
	For the second assertion, let $d=|a-n/2|$.
	Taking $s=\lceil n/2\rceil$ in the hypothesis and using
	$a(n-a)=n^2/4-d^2$, we have
	\[
	\begin{aligned}
		d^2
		&\le Dn+\frac14
		+\ex(a,\cJ(\TP_k))
		-\ex(\lceil n/2\rceil,\cJ(\TP_k))\\
		&\le Dn+\frac14+O_k(dn^{1/3}).
	\end{aligned}
	\]
	The last inequality follows from monotonicity when
	$a\le\lceil n/2\rceil$, and otherwise by summing
	\eqref{eq:increment-bound} from $\lceil n/2\rceil+1$ to $a$.
	Hence $d=O_{k,D}(n^{1/2})$, as required.
\end{proof}

The preceding lemma and \eqref{eq:increment-bound} imply
\[
\max_{0\le s\le n}\{s(n-s)+\ex(s,\cJ(\TP_k))\}
=\left\lfloor\frac{n^2}{4}\right\rfloor
+\ex(\lceil n/2\rceil,\cJ(\TP_k))+O_k(n^{2/3}).
\]
Any maximizing integer $m$ satisfies
$|m-\lceil n/2\rceil|=O_k(n^{1/3})$.
Summing \eqref{eq:increment-bound} between these two integers gives
$ \ex(m,\cJ(\TP_k))
=\ex(\lceil n/2\rceil,\cJ(\TP_k))+O_k(n^{2/3})$.
The upper bound follows from
$m(n-m)\le\lfloor n^2/4\rfloor$,
and the lower bound follows by taking $s=\lceil n/2\rceil$.

To prove \Cref{thm:structure-intro}, we first establish
\Cref{thm:abstract-product} for a class of $3$-chromatic graphs. It shows
that $\ex(n,H)$ exceeds the lower bound in
\Cref{prop:join-criterion-general} by at most $O_H(n)$ and that every
extremal graph differs in at most $O_H(n)$ edges from the complete join of a $\cJ(H)$-free graph and an independent set. Applying this result to $\TP_k$ proves \Cref{thm:structure-intro}.

\begin{theorem}\label{thm:abstract-product}
	Let $H$ be a $3$-chromatic graph.  Suppose that
	\begin{enumerate}[label=\textup{(\roman*)}]
		\item for some $1<\rho<2$ and constants $c,~C>0$,
		$c m^\rho\le\ex(m,\cJ(H))\le C m^\rho$ for all sufficiently large $m$;
		\item there is a partition $V(H)=X\cup Y$ such that both $H[X]$ and
		$H[Y]$ are forests containing at least one edge.
	\end{enumerate}
	Then there exists a constant $C_H>0$, depending only
	on $H,c,C,\rho$, such that, for all sufficiently large $n$,
	\[
	0\le\ex(n,H)-
	\max_{0\le m\le n}\{m(n-m)+\ex(m,\cJ(H))\}
	\le C_Hn.
	\]
	Moreover, every $G\in\EX(n,H)$ has a partition
	$V(G)=A\cup B$, with $a=|A|$, such that $G[A]$ is $\cJ(H)$-free,
	$|a-n/2|\le C_Hn^{\rho/2}$, and
	\[
	e(G[B])+a(n-a)-e_G(A,B)+\ex(a,\cJ(H))-e(G[A])\le C_Hn.
	\]
	There is also a set $Z\subseteq B$, with $|Z|\le C_H$, such that
	$G[B\setminus Z]$ is $H[Y]$-free.  In particular, $e(G[B])+\bigl(|A||B|-e_G(A,B)\bigr)\le C_Hn$.
\end{theorem}

We prove \Cref{thm:abstract-product} using
\Cref{thm:ES-structure,lem:sampling-general}. We first show the properties of the almost balanced bipartition in \Cref{thm:ES-structure} that will be used in the proof of \Cref{thm:abstract-product}.

\begin{lemma}\label{thm:ES-structure}
	Let $H$ be a $3$-chromatic graph.
	There is a nonnegative sequence
	$\zeta_n=\zeta_H(n)\to0$, depending only on $H$, such that,
	for all sufficiently large $n$, every $G\in\EX(n,H)$ admits a partition
	$V(G)=A_1\cup A_2$ minimizing $e(G[A_1])+e(G[A_2])$ and satisfying
	the following properties.
	\begin{enumerate}[label=\textup{(\roman*)}]
		\item For $i=1,2$,
		$\bigl||A_i|-n/2\bigr|\le\zeta_n n$ and $\delta(G)\ge(1/2-\zeta_n)n.$
		
		\item For every $\eta>0$, each $A_i$ contains at most
		$O_{H,\eta}(1)$ vertices $x$ with $d_G(x,A_i)\ge\eta n$.
		\item Let $M$ be a nonempty graph and $t\ge1$ a integer
		such that $H\subseteq M\join I_t$.  For every
		$0<\eta<1/(2v(M))$ and all sufficiently large $n$, let
		\[
		A_i^*(\eta)=\{x\in A_i:d_G(x,A_i)<\eta n\}
		\]
		Then $G[A_i^*(\eta)]$ is $M$-free for $i=1,2$.
	\end{enumerate}
\end{lemma}

\begin{proof}
	
	Parts~\textup{(i)} and~\textup{(ii)} follow from
	Theorem~2.8 of~\cite{GriggsSimonovitsThomas1998}, applied with
	forbidden family $\{H\}$ and $\chi(H)-1=2$.
	
	Minimality of $e(G[A_1])+e(G[A_2])$ implies that moving a vertex
	to the other part cannot decrease this sum.  Together with
	part~\textup{(i)}, this gives, for $i\in\{1,2\}$ and $x\in A_i$,
	\begin{equation}\label{eq:mincut-inequalities}
		\begin{gathered}
			d_G(x,A_{3-i})\ge d_G(x,A_i),\\
			|A_{3-i}\setminus N_G(x)|
			=|A_{3-i}|-d_G(x)+d_G(x,A_i)
			\le d_G(x,A_i)+2\zeta_n n.
		\end{gathered}
	\end{equation}

	To prove~\textup{(iii)}, fix $M,t,\eta$ as stated and suppose there is
	a copy $M'$ of $M$ lies in $G[A_i^*(\eta)]$ for some $i\in\{1,2\}$.
	By \eqref{eq:mincut-inequalities} and the definition of $A_i^*(\eta)$,
	each vertex of $M'$ misses at most $(\eta+2\zeta_n)n$ vertices of
	$A_{3-i}$.  Since $1/2-v(M)\eta>0$ and $\zeta_n\to0$,
	their common neighborhood in $A_{3-i}$ has size at least $ |A_{3-i}|-v(M)(\eta+2\zeta_n)n
	\ge
	\bigl(1/2-v(M)\eta-(2v(M)+1)\zeta_n\bigr)n>t$.  Together with $M'$, any $t$ of these common neighbors
	contain $M\join I_t$ as a subgraph,  hence we find a copy of $H$ in $G$,
	a contradiction.
\end{proof}

For every forest $F$ containing an edge,
we have $\ex(m,F)\le(v(F)-2)m$.
A graph of minimum degree at least $v(F)-1$ contains $F$
by a greedy embedding.
Thus every nonempty subgraph of an $F$-free graph has a vertex
of degree at most $v(F)-2$.
Successively deleting such vertices proves the bound.

We shall combine this bound with the following estimate for the edges
incident with a small vertex set.

\begin{lemma}\label{lem:sampling-general}
	
	Let $\mathcal F$ be a finite family of graphs, and suppose that
	$\ex(m,\mathcal F)\le Cm^\rho$ for every $m\ge1$, where $1<\rho<2$.
	Let $R$ be an $m$-vertex $\mathcal F$-free graph with $m\ge1$,
	and let $U\subseteq V(R)$.
	Define $W=V(R)\setminus U$ and $w=|W|$.
	If $w\le m/2$ and $e(R[U])\le Lm$ for some $L\ge0$, then
	\[
	e(R)\le Lm+C'\left(\frac wm\right)^{\rho-1}m^\rho,
	\]
	where $C'=2^\rho C$.
\end{lemma}

\begin{proof}
	
	If $w=0$, then $U=V(R)$ and $e(R)=e(R[U])\le Lm$,
	as required.
	Suppose $w>0$.  Since $|U|=m-w\ge w$, we can choose
	a uniformly random $w$-subset $U'$ of $U$.
	Every edge of $R[W]$ belongs to $R[W\cup U']$,
	while each edge between $W$ and $U$ belongs to this
	subgraph with probability $w/(m-w)\ge w/m$.
	By linearity of expectation, the expected number of edges
	in $R[W\cup U']$  satisfies
	\[
	\frac wm\bigl(e(R)-e(R[U])\bigr)
	\le \mathbb E\bigl[e(R[W\cup U'])\bigr]
	\le \ex(2w,\mathcal F)
	\le C(2w)^\rho.
	\]
	The upper bound holds because $R[W\cup U']$ is
	$\mathcal F$-free and has $2w$ vertices for every choice
	of $U'$.
	Dividing by $w/m$ and using $e(R[U])\le Lm$, we obtain
	\[
	e(R)
	\le Lm+2^\rho Cmw^{\rho-1}
	=Lm+2^\rho C\left(\frac wm\right)^{\rho-1}m^\rho.
	\]
	This proves the lemma with $C'=2^\rho C$.
\end{proof}

By \Cref{thm:ES-structure,lem:sampling-general}, after deleting
$O_H(1)$ vertices, one part of the partition is $H[Y]$-free. Since
$H[Y]$ is a forest, this part has only $O_H(n)$ edges. Comparing
$G$ with the construction in \Cref{prop:join-criterion-general} then
shows that at most $O_H(n)$ edges between the two parts are missing.

\begin{proof}[Proof of \Cref{thm:abstract-product}]
	
	By increasing $C$ to cover the finitely many smaller values of $m$,
	we may assume that $\ex(m,\cJ(H))\le Cm^\rho$ for every $m\ge1$.

	Let $G$ be an extremal $H$-free graph on $n$ vertices.
	\Cref{prop:join-criterion-general} gives the lower bound in the
	theorem.  Taking $m=\lfloor n/2\rfloor$ and using
	hypothesis~\textup{(i)} also gives $e(G)\ge\left\lfloor n^2/4\right\rfloor+c_*n^\rho$
	for some constant $c_*>0$.
	Let $C'=2^\rho C$ be the constant in
	\Cref{lem:sampling-general}.
	Choose $\eps>0$ such that
	$\eps<1/(32v(H))$ and $C'\bigl(8v(H)\eps\bigr)^{\rho-1}<c_*/4$.
	
	Let $V(G)=A_1\cup A_2$ be the partition given by
	\Cref{thm:ES-structure}.
	We remove from each part the vertices with at least $\eps n$
	neighbors in that part.  Write $ Z_i=\{x\in A_i:d_G(x,A_i)\ge\eps n\}$ and
	$A_i^*=A_i\setminus Z_i$, where $i=1,2$. By \Cref{thm:ES-structure}\textup{(ii)}, $|Z_1|+|Z_2|=O_H(1)$, so only boundedly many vertices are removed.
	Every $M\in\cJ(H)$ satisfies
	$v(M)\le v(H)$ and $H\subseteq M\join I_{v(H)}$.
	Our choice of $\eps$ ensures that $\eps<1/(2v(M))$
	for every such $M$.
	Applying \Cref{thm:ES-structure}\textup{(iii)} with
	$\eta=\eps$ to each member of the finite family $\cJ(H)$
	therefore shows that both $G[A_1^*]$ and $G[A_2^*]$
	are $\cJ(H)$-free.
	
	Since $e_G(A_1,A_2)\le\lfloor n^2/4\rfloor$,
	the lower bound on $e(G)$ implies
	$e(G[A_1])+e(G[A_2])\ge c_*n^\rho$.
	Deleting $Z_1\cup Z_2$ removes only $O_H(n)$ edges, so
	$e(G[A_1^*])+e(G[A_2^*])\ge c_*n^\rho-O_H(n)$. Choose the labels such that $e(G[A_1^*])\ge e(G[A_2^*])$.
	Since $\rho>1$, the preceding estimate yields
	$e(G[A_1^*])
	\ge \frac{c_*}{2}n^\rho-O_H(n)
	\ge \frac{c_*}{3}n^\rho$.

	Next, we show $G[A_2^*]$ is $H[Y]$-free.
	Suppose not, that $G[A_2^*]$ contains a copy of
	$H[Y]$ with vertex set $Y'$.  Let
	$U=A_1^*\cap\bigcap_{y\in Y'}N_G(y)$,
	$W=A_1^*\setminus U$,
	and $m_1=|A_1^*|$.
	Each $y\in Y'$ has fewer than $\eps n$ neighbors in $A_2$,
	so \eqref{eq:mincut-inequalities} gives
	$|W|\le v(H[Y])(\eps+2\zeta_H(n))n$.
	Since $m_1=n/2+o(n)$ and $\zeta_H(n)\to0$, we obtain
	$|W|/m_1\le8v(H)\eps<1/2$.
	
	Since $G$ is $H$-free and $U$ is completely joined to $Y'$, $G[U]$ is $H[X]$-free.
	Since $H[X]$ is a forest,
	$e(G[U])\le(v(H[X])-2)|U|\le(v(H[X])-2)m_1$.
	Applying \Cref{lem:sampling-general} to $G[A_1^*]$ yields
	\[
	\begin{aligned}
		e(G[A_1^*])
		&\le(v(H[X])-2)m_1
		+C'\bigl(8v(H)\eps\bigr)^{\rho-1}m_1^\rho\\
		&\le(v(H)-2)n+\frac{c_*}{4}n^\rho
		<\frac{c_*}{3}n^\rho,
	\end{aligned}
	\]
	contradicting the lower bound on $e(G[A_1^*])$.
	Thus $G[A_2^*]$ is $H[Y]$-free, and
	$e(G[A_2^*])\le(v(H[Y])-2)|A_2^*|=O(n)$.
	
	Let $A=A_1^*$, $B=V(G)\setminus A$, $a=|A|$, and
	$Z=Z_1\cup Z_2$. Then $G[A]$ is $\cJ(H)$-free, $Z\subseteq B$, $|Z|=O_H(1)$, and $B\setminus Z=A_2^*$ induces an $H[Y]$-free graph.
	Since $B=A_2^*\cup Z$, we also have $e(G[B])=O_H(n)$.
	
	Since $G\in\EX(n,H)$, \Cref{prop:join-criterion-general} gives 
	$e(G)\ge a(n-a)+\ex(a,\cJ(H))$.
	Consequently,
	\[
	0\le
	a(n-a)-e_G(A,B)+\ex(a,\cJ(H))-e(G[A])
	\le e(G[B])=O(n).
	\]
	
	Adding $e(G[B])$ gives the inequality claimed for the partition
	$A\cup B$.  Since $e(G[B])=O(n)$ and $G[A]$ is $\cJ(H)$-free,
	\[
	\begin{aligned}
		\ex(n,H)=e(G)
		&\le a(n-a)+\ex(a,\cJ(H))+O_H(n)\\
		&\le
		\max_{0\le m\le n}
		\{m(n-m)+\ex(m,\cJ(H))\}+O_H(n).
	\end{aligned}
	\]
	Using $e(G)\ge\lfloor n^2/4\rfloor$ in the first inequality,
	we get $ \left(a-\frac n2\right)^2
	\le \ex(a,\cJ(H))+O_H(n)
	\le Cn^\rho+O_H(n)$.
	Therefore $|a-n/2|=O(n^{\rho/2})$ holds.
	
	Deleting the edges within $B$ and adding the missing edges between
	$A$ and $B$ transforms $G$ into $G[A]\join I_{|B|}$ using
	$e(G[B])+a(n-a)-e_G(A,B)=O_H(n)$ edge modifications.
	Taking $C_H$ sufficiently large then completes the proof.
\end{proof}

\medskip
The linear error in \Cref{thm:abstract-product} cannot, in general,
be replaced by $o(n)$. Take $H=K_{2,2,2}$.
Every independent set of $H$ is contained in a vertex class,
whose deletion produces a $C_4$; hence $\cJ(H)=\{C_4\}$.
Moreover, $\ex(m,C_4)=\Theta(m^{3/2})$, and there is a partition $V(H)=X\cup Y$ such that $H[X]\cong H[Y]\cong P_3$. Hence $H$ satisfies the two assumptions of \Cref{thm:abstract-product}.

By Theorem~2.6 of~\cite{GriggsSimonovitsThomas1998}, $\ex(n,K_{2,2,2})
=
\max_{0\le m\le n}
\left\{
m(n-m)+\ex(m,C_4)
+\left\lfloor (n-m)/2\right\rfloor
\right\}$ holds
for all sufficiently large $n$.
Let $m$ maximize $m(n-m)+\ex(m,C_4)$.
Since
$m(n-m)+\ex(m,C_4)\ge\lfloor n^2/4\rfloor$
and $\ex(m,C_4)=O(n^{3/2})$, we have
$m=n/2+O(n^{3/4})$. Then $\ex(n,K_{2,2,2})
-
\max_{0\le m\le n}\{m(n-m)+\ex(m,C_4)\}
\ge \frac n4-O(n^{3/4})$.
Therefore, the upper bound in \Cref{thm:abstract-product} shows
that this difference is $\Theta(n)$.

We prove \Cref{thm:structure-intro} by applying
\Cref{thm:abstract-product} to $\TP_k$.
The second assertion of \Cref{lem:join-maximizers} gives
the sharper bound on the part sizes.

\begin{proof}[Proof of \Cref{thm:structure-intro}]
	
	Partition $V(\TP_k)$ according to the parity of the layer index.
	Each part induces a disjoint union of paths, since edges between
	distinct layers join consecutive layers.
	The even and odd parts contain $P_2$ and $P_3$, respectively,
	so both forests contain an edge.
	Together with \Cref{prop:join-family-comparison}, this verifies
	the hypotheses of \Cref{thm:abstract-product} with $\rho=4/3$.
	
	\Cref{thm:abstract-product} gives the required two-sided estimate.
	It also gives, for every $G\in\EX(n,\TP_k)$, a partition
	$V(G)=A\cup B$ such that $G[A]$ is $\cJ(\TP_k)$-free,
	$e(G[B])=O_k(n)$, and $G$ differs from
	$G[A]\join I_{|B|}$ in $O_k(n)$ edges.
	Let $a=|A|$. The edge count and
	\Cref{prop:join-criterion-general} give
	\[
	\max_{0\le m\le n}
	\{m(n-m)+\ex(m,\cJ(\TP_k))\}
	\le e(G)
	\le a(n-a)+\ex(a,\cJ(\TP_k))+O_k(n).
	\]
	The second assertion of \Cref{lem:join-maximizers} therefore gives
	$|a-n/2|=O_k(n^{1/2})$.
	Taking $R=G[A]$ and $m=a$ proves the structural conclusion.
\end{proof}

The estimate following \Cref{lem:join-maximizers}, together with
\Cref{thm:structure-intro}, shows
\[
\ex(n,\TP_k)
=\left\lfloor\frac{n^2}{4}\right\rfloor
+\ex(\lceil n/2\rceil,\cJ(\TP_k))+O_k(n).
\]
In particular, \Cref{prop:join-family-comparison} gives
$\ex(n,\TP_k)=\lfloor n^2/4\rfloor+\Theta_k(n^{4/3})$.

\section{\texorpdfstring{The exact result for $\TP_4$}{The exact result for TP4}}
\label{sec:exact-four-layer}

In this section, we will prove Theorem~\ref{exactp4}.
We determine $\cJ(\TP_4)$, establish a density gap for
$\cJ(\TP_4)$-free graphs, and finally show that one part of the
structural partition is independent.

\begin{figure}[htbp]
	\centering
	\begin{tikzpicture}[scale=.61,every node/.style={circle,inner sep=1.8pt,fill=black}]
		\begin{scope}[xshift=-5.4cm]
			\foreach \i in {1,...,5}{\foreach \j in {1,...,\i}{
					\coordinate (a\i\j) at
					({1.05*(\j-(\i+1)/2)},{-1.02*(\i-1)});
			}}
			\draw (a21)--(a22)--(a33)--(a43)--(a42)--(a31)--(a21);
			\draw (a42)--(a52)--(a51);
			\draw (a43)--(a54)--(a55);
			\foreach \v in {21,22,31,33,42,43,51,52,54,55}{\node at (a\v) {};}
			\foreach \v in {21,22,31,33,42,43,51,52,54,55}{
				\pgfmathtruncatemacro{\rr}{int(\v/10)}
				\pgfmathtruncatemacro{\ss}{mod(\v,10)}
				\node[fill=none,draw=none,rectangle,inner sep=0pt,font=\tiny]
				at ($(a\v)+(0,-.50)$) {$v_{\ss}^{\rr}$};
			}
			\node[fill=none,draw=none,rectangle] at (0,-4.55) {$F_1$};
		\end{scope}
		\begin{scope}
			\foreach \i in {1,...,5}{\foreach \j in {1,...,\i}{
					\coordinate (b\i\j) at ({1.05*(\j-(\i+1)/2)},{-0.88*(\i-1)});
			}}
			\draw (b11)--(b21)--(b32)--(b33)--(b43)--(b54)--(b55);
			\draw (b32)--(b42)--(b41)--(b52)--(b42);
			\draw (b32)--(b43)--(b42)--(b32);
			\foreach \v in {11,21,32,33,41,42,43,52,54,55}{\node at (b\v) {};}
			\foreach \v in {11,21,32,33,41,42,43,52,54,55}{
				\pgfmathtruncatemacro{\rr}{int(\v/10)}
				\pgfmathtruncatemacro{\ss}{mod(\v,10)}
				\node[fill=none,draw=none,rectangle,inner sep=0pt,font=\tiny]
				at ($(b\v)+(0,-.50)$) {$v_{\ss}^{\rr}$};
			}
			\node[fill=none,draw=none,rectangle] at (0,-4.55) {$F_2$};
		\end{scope}
		\begin{scope}[xshift=5.4cm]
			\foreach \i in {1,...,5}{\foreach \j in {1,...,\i}{
					\coordinate (c\i\j) at ({1.05*(\j-(\i+1)/2)},{-0.88*(\i-1)});
			}}
			\draw (c21)--(c22)--(c32)--(c21);
			\draw (c32)--(c42)--(c43)--(c32);
			\draw (c41)--(c42)--(c52)--(c41);
			\draw (c43)--(c44)--(c54)--(c43);
			\foreach \v in {21,22,32,41,42,43,44,52,54}{\node at (c\v) {};}
			\foreach \v in {21,22,32,41,42,43,44,52,54}{
				\pgfmathtruncatemacro{\rr}{int(\v/10)}
				\pgfmathtruncatemacro{\ss}{mod(\v,10)}
				\node[fill=none,draw=none,rectangle,inner sep=0pt,font=\tiny]
				at ($(c\v)+(0,-.50)$) {$v_{\ss}^{\rr}$};
			}
			\node[fill=none,draw=none,rectangle] at (0,-4.55) {$F_3$};
		\end{scope}
	\end{tikzpicture}
	\caption{The graphs $F_1,F_2,F_3$.}
	\label{fig:obstructions}
\end{figure}
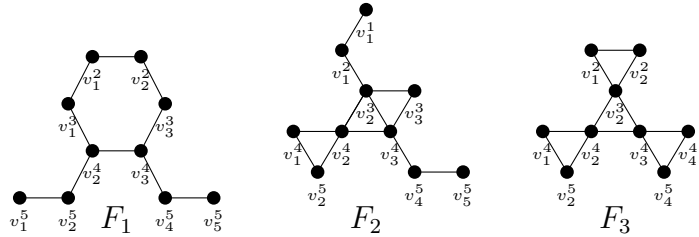

Define
\begin{align*}
	I_1=\{v_1^1,v_2^3,v_1^4,v_4^4,v_3^5\},~I_2=\{v_2^2,v_1^3,v_4^4,v_1^5,v_3^5\},~\text{and}~
	I_3=\{v_1^1,v_1^3,v_3^3,v_1^5,v_3^5,v_5^5\}.
\end{align*}
By the definition of $\TP_k$, we have each
$I_i$ is independent.  Let $F_i=\TP_4- I_i$ for
$i=1,2,3$.  The three graphs are shown in \Cref{fig:obstructions}.

Since $I_1=V_2$, we have
$F_1=B_{4,2}$.  We now determine $\cJ(\TP_4)$. This is the only
computer-assisted part of the proof and is verified by the computation in
\Cref{app:verification}.

\begin{proposition}\label{prop:classification}
	$\cJ(\TP_4)=\{F_1,F_2,F_3\}$.
\end{proposition}

\begin{proof}
	It suffices to consider maximal independent sets: if $I$ is not maximal,
	extend it to a maximal independent set $I'$, and then
	$\TP_4-I'\subseteq\TP_4-I$.
	
	The computation finds $51$ maximal independent sets, which form
	$12$ symmetry classes. For one representative of each class, the program
	verifies that $\TP_4-I$ contains at least one of $F_1,F_2,F_3$; the
	representatives are listed in
	\Cref{tab:maximal-independent} in the appendix.
	
	It remains to show that these three graphs are pairwise incomparable.
	The graph $F_1$ is bipartite, whereas $F_2$ and $F_3$ contain triangles,
	so neither is a subgraph of $F_1$. Since $F_1$ has ten vertices and
	contains a copy of $C_6$, it is a subgraph of neither the nine-vertex
	graph $F_3$ nor the $C_6$-free graph $F_2$. Finally, $F_2$ has more
	vertices than $F_3$, while $F_3$ contains four triangles and $F_2$
	contains only three. Thus neither $F_2$ nor $F_3$ is a subgraph of the
	other.
	
	Therefore, every graph $\TP_4-I$ contains one of $F_1,F_2,F_3$, and
	these three graphs are precisely the subgraph-minimal members. Hence $\cJ(\TP_4)=\{F_1,F_2,F_3\}$.
\end{proof}

We next show that an $\{F_1,F_2,F_3\}$-free graph containing a vertex of degree at least $\eta m$ has at most $\ex(m,\{F_1,F_2,F_3\})-\gamma m^{4/3}$ edges. This estimate will later be used to show that the exceptional set in the structural partition is empty.

\begin{lemma}\label{lem:density-gap}
	For every $\eta>0$, there are constants $\gamma=\gamma(\eta)>0$ and
	$m_0=m_0(\eta)$ such that the following holds.  If $R$ is a $\{F_1,F_2,F_3\}$-free graph
	on $m\ge m_0$ vertices and some $x\in V(R)$ has $d_R(x)\ge\eta m$, then
	$e(R)\le \ex(m,\{F_1,F_2,F_3\})-\gamma m^{4/3}$.
\end{lemma}
\begin{proof}
	The graph $F_1$ has the unique cycle $v_1^2v_2^2v_3^3v_3^4v_2^4v_1^3v_1^2$.
	Let $S=F_1-v_3^3$. Then $S$ is a tree on nine vertices and hence is bipartite. Since $N_{F_1}(v_3^3)=\{v_2^2,v_3^4\}$, and $v_2^2$ and $v_3^4$ are at distance four in $S$, they belong to the same part. Denote this part by $P$.
	
	Let $N=N_R(x)$ and $W=V(R)\setminus\bigl(N\cup\{x\}\bigr)$. The graph $R[N]$ contains no copy of $S$, since such a copy together with $x$ would contain $F_1$. Because $S$ is a tree on nine vertices, every graph of minimum degree at least eight contains a copy of $S$ by a greedy construction. Consequently, every subgraph of $R[N]$ has a
	vertex of degree at most seven, and hence $e(R[N])\le 7|N|\le 7m$. 
	
	We now bound the number of edges between $N$ and $W$. Repeatedly delete a vertex of $N$ having at most three neighbors in $W$, or a vertex of $W$ having at most four neighbors in $N$. This process must delete every vertex. Otherwise, the remaining graph has minimum degree at least four from $N$ to $W$ and at least five from $W$ to $N$, and we can construct a copy of $S$ with $P$ in $N$, a contradiction. Since every edge between $N$ and $W$ is counted when its first endpoint is deleted, there are at most $3|N|+4|W|=O(m)$ such edges.
	
	The graph $R[W]$ is $\{F_1,F_2,F_3\}$-free, and $x$ has no neighbors in $W$. Hence, for some constant $C>0$, $e(R)\le \ex\bigl(|W|,\{F_1,F_2,F_3\}\bigr)+Cm$. 
	
	By \Cref{prop:classification,prop:join-family-comparison}, there is a constant $c_0>0$ such that $\ex\bigl(u,\{F_1,F_2,F_3\}\bigr)\ge c_0u^{4/3}$ for all sufficiently large $u$. Moreover, $m-|W|=|N|+1\ge\eta m$. Since $F_1,F_2,F_3$ are connected, the disjoint union of two $\{F_1,F_2,F_3\}$-free graphs is also $\{F_1,F_2,F_3\}$-free. Therefore, \[ \begin{aligned} \ex\bigl(m,\{F_1,F_2,F_3\}\bigr) &\ge \ex\bigl(|W|,\{F_1,F_2,F_3\}\bigr) +\ex\bigl(m-|W|,\{F_1,F_2,F_3\}\bigr)\\ &\ge \ex\bigl(|W|,\{F_1,F_2,F_3\}\bigr) +c_0\eta^{4/3}m^{4/3}. \end{aligned} \] Consequently, \[ e(R)\le \ex\bigl(m,\{F_1,F_2,F_3\}\bigr) -c_0\eta^{4/3}m^{4/3}+Cm. \]
	Let $\gamma=\frac12c_0\eta^{4/3}$. This completes the proof of Lemma~\ref{lem:density-gap}.
\end{proof}

Now we will prove~\Cref{exactp4}.

We need the following decomposition of $\TP_4$ in the proof of
\Cref{exactp4}. Let
\[
T=\TP_4-\{v_2^2,v_1^3,v_3^4,v_2^5,v_4^5\}.
\]
Then $T$ is a tree on ten vertices, and the only edge among the five
deleted vertices is $v_3^4v_4^5$. In the proof, this decomposition will
be used to show that the common neighborhood of a possible edge contains
no copy of $T$.

\begin{proof}[Proof of \Cref{exactp4}]
	
	Let $G\in\EX(n,\TP_4)$.
	By \Cref{prop:classification,prop:join-criterion-general,prop:join-family-comparison},
	there is a constant $c_*>0$ such that
	\[
	\begin{aligned}
		\left\lfloor\frac{n^2}{4}\right\rfloor+c_*n^{4/3}
		&\le
		\max_{0\le m\le n}
		\{m(n-m)+\ex(m,\{F_1,F_2,F_3\})\}\le e(G).
	\end{aligned}
	\]
	Let $C>0$ satisfy
	$\ex(m,\{F_1,F_2,F_3\})\le Cm^{4/3}$ for every $m\ge1$,
	and let $C'=2^{4/3}C$ as in \Cref{lem:sampling-general}.
	Choose $\eps>0$ satisfying
	$\eps<1/100$, and $C'(6\eps)^{1/3}< c_*/10$. Applying \Cref{thm:ES-structure} to $G$ to obtain a partition
	$V(G)=A_1\cup A_2$, and write $\zeta_n=\zeta_{\TP_4}(n)$.
	Let $Z_i=\{x\in A_i:d_G(x,A_i)\ge\eps n\}$ and $A_i^*=A_i\setminus Z_i$, where $i=1,2$.
	
	By \Cref{thm:ES-structure}\textup{(i) and (ii)},
	$|A_i|=n/2+o(n)$ and $|Z_1|+|Z_2|=O(1)$.
	Since $v(F_j)\le10$ and $\TP_4\subseteq F_j\join I_6$
	for $j=1,2,3$, \Cref{thm:ES-structure}\textup{(iii)}
	shows that both $G[A_1^*]$ and $G[A_2^*]$
	are $\{F_1,F_2,F_3\}$-free.
	
	The lower bound on $e(G)$ and
	$e_G(A_1,A_2)\le\lfloor n^2/4\rfloor$ give
	$e(G[A_1])+e(G[A_2])\ge c_*n^{4/3}$.
	Deleting $Z_1\cup Z_2$ removes only $O(n)$ edges. Without loss of generality, suppose $e(G[A_1^*])\ge e(G[A_2^*])$. Then
	$e(G[A_1^*])
	\ge \frac{c_*}{2}n^{4/3}-O(n)
	\ge \frac{c_*}{3}n^{4/3}$.
	We first show that $A_2^*$ is independent.
	Suppose not, there exists $xy\in E(G[A_2^*])$. Let
	$ U=N_G(x)\cap N_G(y)\cap A_1^*,~W=A_1^*\setminus U, ~\text{and}~m_1=|A_1^*|$.
	By the definition of $A_2^*$ and \eqref{eq:mincut-inequalities},
	each of $x,y$ has at most $(\eps+2\zeta_n)n$ non-neighbors
	in $A_1^*$.
	Every vertex of $W$ is nonadjacent to at least one of $x,y$,
	so $|W|\le2(\eps+2\zeta_n)n$.
	Since $m_1=n/2+o(n)$ and $\zeta_n\to0$,
	we obtain $|W|/m_1\le6\eps<1/2$.
	
	We claim that $G[U]$ is $T$-free.
	Suppose not, that $G[U]$ contains a copy $T'$ of $T$.
	Each vertex of $T'$ lies in $A_1^*$, so
	\eqref{eq:mincut-inequalities} gives
	\[
	\left|
	\left(A_2^*\setminus\{x,y\}\right)
	\cap\bigcap_{v\in V(T')}N_G(v)
	\right|
	\ge |A_2^*|-2-10(\eps+2\zeta_n)n>3.
	\]
	The copy $T'$, together with the edge $xy$ and three further vertices from this common neighborhood, forms a copy of $\TP_4$, a contradiction. Since $T$ is a tree on ten vertices, every subgraph of $G[U]$ has a vertex of degree at most eight. Therefore, $e(G[U])\le 8|U|\le 8|A_1^*|=8m_1$ holds. 
	
	Applying \Cref{lem:sampling-general} to $G[A_1^*]$, we get
	\[
	\begin{aligned}
		e(G[A_1^*])
		\le 8m_1+C'(6\eps)^{1/3}m_1^{4/3}
		\le 8n+\frac{c_*}{10}n^{4/3}
		<\frac{c_*}{3}n^{4/3},
	\end{aligned}
	\]
	contrary to the lower bound on $e(G[A_1^*])$.
	Thus $A_2^*$ is independent.
	
	We next use \Cref{lem:density-gap} to show that
	$Z_1\cup Z_2=\varnothing$.
	Suppose that $Z_i\ne\varnothing$ for some $i\in\{1,2\}$.
	Let $x\in Z_i$ and $D=d_G(x,A_i)\ge\eps n$.
	Let
	$ Z'=(Z_1\cup Z_2)\cup
	\bigcup_{j=1}^2
	\{v\in A_j:d_G(v,A_j)\ge D/1000\}$, and $W_j=A_j\setminus Z'
	\quad(j=1,2)$.
	Since $D/1000\ge(\eps/1000)n$,
	\Cref{thm:ES-structure}\textup{(ii)} gives $|Z'|=O(1)$.
	Moreover, $W_j\subseteq A_j^*$, so $G[W_1]$ is
	$\{F_1,F_2,F_3\}$-free and $W_2$ is independent.
	By \eqref{eq:mincut-inequalities}, we have
	$d_G(x,A_{3-i})\ge d_G(x,A_i)=D$.
	It follows that $|N_G(x)\cap W_j|\ge D-|Z'|~\text{for}~j=1,2$.
	
	Let $X_0=W_1\cup\{x\}$.
	We claim that $G[X_0]$ is $\{F_1,F_2,F_3\}$-free.
	Any copy $F'$ of some $F_j$ in $G[X_0]$ must contain $x$,
	since $G[W_1]$ is $\{F_1,F_2,F_3\}$-free.
	Every $y\in V(F')\setminus\{x\}$ lies in $W_1$, and
	\eqref{eq:mincut-inequalities} gives $
	|W_2\setminus N_G(y)|\le D/1000+2\zeta_n n$.
	Using $v(F_j)\le10$, $|N_G(x)\cap W_2|\ge D-|Z'|$,
	and $\zeta_n n=o(D)$, we obtain
	\[
	\left|W_2\cap\bigcap_{v\in V(F')}N_G(v)\right|
	\ge D-|Z'|-9(D/1000+2\zeta_n n)>6.
	\]
	Since $|I_j|\le6$, choose $|I_j|$ vertices from $W_2\cap\bigcap_{v\in V(F')}N_G(v)$. Together with $F'$, these vertices contain a copy of $\TP_4$, a contradiction. This proves the claim.
	
	Now $|X_0|=n/2+o(n)$ and $ d_{G[X_0]}(x)\ge D-|Z'|\ge(\eps/2)|X_0|$.
	By \Cref{lem:density-gap}, there is a constant $\gamma_0>0$
	such that
	$e(G[X_0])
	\le \ex(|X_0|,\{F_1,F_2,F_3\})
	-\gamma_0|X_0|^{4/3}$.
	Let
	$Y_0=V(G)\setminus X_0=W_2\cup(Z'\setminus\{x\})$.
	Since $W_2$ is independent and $|Z'|=O(1)$,
	we have $e(G[Y_0])=O(n)$.
	Consequently,
	\[
	\begin{aligned}
		e(G)
		&\le |X_0||Y_0|+\ex(|X_0|,\{F_1,F_2,F_3\})
		-\gamma_0|X_0|^{4/3}+O(n)\\
		&\le \max_{0\le m\le n}
		\{m(n-m)+\ex(m,\{F_1,F_2,F_3\})\}
		-\gamma_0|X_0|^{4/3}+O(n).
	\end{aligned}
	\]
	Because $|X_0|=n/2+o(n)$, this contradicts the lower bound
	at the start of the proof.
	Therefore $Z_1\cup Z_2=\varnothing$.
	
	It follows that $G[A_1]$ is $\{F_1,F_2,F_3\}$-free
	and $A_2$ is independent.
	Let $a=|A_1|$. By
	\Cref{prop:classification,prop:join-criterion-general}, we get
	\[
	\begin{aligned}
		\ex(n,\TP_4)=e(G)
		&\le a(n-a)+e(G[A_1])\\
		&\le a(n-a)+\ex(a,\{F_1,F_2,F_3\})\\
		&\le \max_{0\le m\le n}
		\{m(n-m)+\ex(m,\{F_1,F_2,F_3\})\}\\
		&\le \ex(n,\TP_4).
	\end{aligned}
	\]
	Thus $G=G[A_1]\join I_{n-a}$,
	$G[A_1]\in\EX(a,\{F_1,F_2,F_3\})$, and the maximum is attained when $m=a$.
	Conversely, if the maximum occurs at $m$, and let
	$R\in\EX(m,\{F_1,F_2,F_3\})$. 
	Then $R\join I_{n-m}$ is $\TP_4$-free and has
	$m(n-m)+\ex(m,\cJ(\TP_4))=\ex(n,\TP_4)$
	edges. Hence it belongs to $\EX(n,\TP_4)$. This completes the proof of~\Cref{exactp4}.
\end{proof}

\section{Concluding remarks and open problems}\label{sec:conclusion}

Our analysis of the subgraphs $B_{k,r}$ establishes
$\ex(n,\cJ(\TP_k))=\Theta_k(n^{4/3})$ for every $k\ge4$.
Combining the partition in \Cref{thm:ES-structure} with
the sampling estimate in \Cref{lem:sampling-general}
and the linear Tur\'an bound for forests, we get
\[
\ex(n,\TP_k)=
\max_{0\le m\le n}
\{m(n-m)+\ex(m,\cJ(\TP_k))\}+O_k(n).
\]
The same argument shows that every sufficiently large extremal
$\TP_k$-free graph differs in $O_k(n)$ edges from
$R\join I_{n-m}$, where $R$ is an $m$-vertex
$\cJ(\TP_k)$-free graph and $m=n/2+O_k(n^{1/2})$.
For $k=4$, the classification of $\cJ(\TP_4)$ and the
estimate in \Cref{lem:density-gap} allow us to eliminate
the exceptional vertices and obtain the exact formula
and the characterization of all extremal graphs in
\Cref{exactp4}.

The order of $\ex(m,\{F_1,F_2,F_3\})$ is determined,
but an asymptotic formula is not known.
Taking a maximum bipartite subgraph of an extremal
$C_6$-free graph gives the lower bound below.
The upper bound follows from the structure of $F_1$
by deleting low-degree vertices and extending a copy of
$C_6$, as in \cite{ChenChenZhu2026}.
We get
\[
\frac12\ex(m,C_6)
\le \ex(m,\{F_1,F_2,F_3\})
\le \ex(m,C_6)+O(m).
\]

We conclude with several questions arising from these results,
concerning sharper Tur\'an estimates and the structure of
extremal graphs.

\begin{problem}\label{prob:four-layer-family}
	Does the limit $\lim_{m\to\infty}\frac{
		\ex(m,\{F_1,F_2,F_3\})}{m^{4/3}}$
	exist, and if so, what is its value?
	Is $ \ex(m,\{F_1,F_2,F_3\})
	=(1+o(1))\ex(m,C_6)?$
\end{problem}

For $k\ge5$, \Cref{thm:structure-intro} leaves a
nonnegative error of order at most $n$.

\begin{problem}\label{prob:linear-correction}
	For each integer $k\ge5$, determine whether
	\[
	\ex(n,\TP_k)=
	\max_{0\le m\le n}
	\{m(n-m)+\ex(m,\cJ(\TP_k))\}
	\]
	holds for all sufficiently large $n$. If this equality fails, determine the
	order of the difference between the two sides.
\end{problem}

The next question asks for all extremal graphs.
An affirmative answer would also establish the equality
in \Cref{prob:linear-correction}.

\begin{problem}\label{prob:exact-extremal}
	For $k\ge5$, is every sufficiently large extremal
	$\TP_k$-free graph of the form $R\join I_{n-m}$, where
	$R\in\EX(m,\cJ(\TP_k))$ and $m$ maximizes
	$m(n-m)+\ex(m,\cJ(\TP_k))$?
	If not, characterize the remaining extremal graphs.
\end{problem}

The proof for $k=4$ uses \Cref{lem:density-gap}
to exclude vertices of large internal degree.
The same estimate holds for $\cJ(\TP_5)$:
every $\cJ(\TP_5)$-free graph is $B_{5,2}$-free,
and $B_{5,2}-v_2^4$ is a tree on thirteen vertices.
The neighborhood estimates and the superadditivity argument
in the proof of \Cref{lem:density-gap} therefore apply
with $B_{5,2}$ in place of $F_1$.
This leads to the following question for larger $k$.

\begin{problem}\label{prob:density-gap}
	For every  $k\ge6$ and $0<\eta<1$, does there exist
	constants $\gamma>0$ and $m_0=m_0(k,\eta)$, such that every $m$-vertex $\cJ(\TP_k)$-free
	graph $R$ with $m\ge m_0$ and $\Delta(R)\ge\eta m$ satisfies
	$e(R)\le\ex(m,\cJ(\TP_k))-\gamma m^{4/3}?$
\end{problem}

\appendix
\section{Finite verification and reproducibility}\label{app:verification}

\Cref{prop:classification} is the only
computer-assisted ingredient in the proof of \Cref{exactp4}.
The script \texttt{verification/verify\_tp4.py} constructs $\TP_4$ from
Definition~\ref{TPk}, enumerates all $2^{15}$ vertex subsets, and performs,
in particular, the following checks.
\begin{enumerate}[label=\textup{(\roman*)}]
	
	\item It finds the $51$ maximal independent sets and the twelve symmetry
	classes represented in \Cref{tab:maximal-independent}.
	
	\item It verifies that every listed $\TP_4-I$ contains the stated members of
	$\{F_1,F_2,F_3\}$.
	\item It verifies all required noncontainments, including the pairwise
	incomparability of $F_1,F_2,F_3$.
\end{enumerate}
The program is deterministic and takes no external
input.  Its transcript, the program, a README, and a SHA-256 manifest are
included in the supplementary archive.

For completeness, the following table lists one representative from
each symmetry class of maximal independent sets of $\TP_4$ and records
which of $F_1,F_2,F_3$ occur in $\TP_4-I$.

\begin{table}[H]
	\centering
	\scriptsize
	\setlength{\tabcolsep}{2pt}
	\renewcommand{\arraystretch}{0.95}
	\begin{tabular}{|
			>{\centering\arraybackslash}p{.46\textwidth}|
			>{\centering\arraybackslash}p{.12\textwidth}|
			>{\centering\arraybackslash}p{.34\textwidth}|}
		\hline
		Representative $I$ & Class size & Contained graphs\\ \hline
		$\{v_1^2,v_4^4,v_2^5\}$ & $2$ & $F_1,F_2,F_3$\\ \hline
		$\{v_1^1,v_1^3,v_4^4,v_2^5\}$ & $6$ & $F_1,F_2$\\ \hline
		$\{v_1^1,v_2^3,v_1^4,v_4^5\}$ & $6$ & $F_1$\\ \hline
		$\{v_1^1,v_2^3,v_1^5,v_4^5\}$ & $6$ & $F_1$\\ \hline
		$\{v_1^1,v_2^3,v_2^5,v_4^5\}$ & $3$ & $F_1$\\ \hline
		$\{v_1^2,v_3^3,v_1^4,v_4^5\}$ & $6$ & $F_1$\\ \hline
		$\{v_1^1,v_1^3,v_3^3,v_1^5,v_4^5\}$ & $6$ & $F_1,F_2$\\ \hline
		$I_2=\{v_2^2,v_1^3,v_4^4,v_1^5,v_3^5\}$ & $3$ & $F_2$\\ \hline
		$\{v_1^1,v_1^3,v_3^4,v_1^5,v_5^5\}$ & $3$ & $F_1$\\ \hline
		$\{v_1^1,v_1^3,v_3^4,v_2^5,v_5^5\}$ & $6$ & $F_1$\\ \hline
		$I_1=\{v_1^1,v_2^3,v_1^4,v_4^4,v_3^5\}$ & $3$ & $F_1$\\ \hline
		$I_3=\{v_1^1,v_1^3,v_3^3,v_1^5,v_3^5,v_5^5\}$ & $1$ & $F_3$\\ \hline
	\end{tabular}
	\caption{Maximal independent sets of $\TP_4$ up to symmetry.}
	\label{tab:maximal-independent}
\end{table}

\section*{Code availability}
The verification script and its output transcript are
provided in the supplementary material.
The script requires Python~3.10 or later and uses only
the standard library. No external data are required.

\section*{Acknowledgments}

\noindent {\bf AI-assisted technologies.} ChatGPT 6 were used to improve the presentation of the manuscript. All AI-generated suggestions were verified and refined by the authors, who take full responsibility for the correctness and originality of the paper.\\
 
\noindent {\bf Funding.}  Xiao was  founded by  Basic Research Program of Jiangsu (Grants No. BK20250673).  Yuan was founded by National Natural Science Foundation of China (No. 12271169 and 12331014) and Science and Technology Commission of Shanghai Municipality (No. 22DZ2229014).

\small
\setstretch{1.0}

\end{document}